\documentclass[a4paper, 12pt]{article}
\usepackage{amsfonts}
\usepackage{amssymb}
\usepackage{amsmath}
\usepackage{eufrak}
\usepackage{latexsym}
\usepackage{amsbsy}
\usepackage[all]{xy}

\usepackage[latin1]{inputenc}
\usepackage[frenchb]{babel}
\usepackage[T1]{fontenc}     

\newtheorem{defn}{D\'efinition}[section]
\newtheorem{thm}[defn]{Th\'eor\`eme}
\newtheorem{prop}[defn]{Proposition}
\newtheorem{lem}[defn]{Lemme}

\newtheorem{rem}[defn]{Remarque}

\numberwithin{equation}{section}

\newcommand{\proof}{\noindent {\textsc{Démonstration :} }\rm}
\newcommand{\fin}{\hfill{\Large$\Box$}\par}

\newcommand{\al}{\alpha}

\newcommand{\vphi}{\varphi}

\newcommand{\la}{\lambda}
\newcommand{\Ga}{\Gamma}
\newcommand{\ga}{\gamma}
\newcommand{\si}{\sigma}
\newcommand{\Om}{\Omega}
\newcommand{\om}{\omega}

\newcommand{\Sg}{\mathbb {S}}
\newcommand{\C}{\mathbb {C}}

\newcommand{\Z}{\mathbb {Z}}

\newcommand{\U}{\mathbb {U}}
\newcommand{\Pj}{\mathbb {P}}

\def\abs#1{\vert #1\vert}
\def\bignorm#1{\left\|\, #1\,\right\|}
\def\bigabs#1{\left\vert\, #1\,\right\vert}

\def\CC{{\cal C}}

\def\DD{{\cal D}}

\def\F{{\cal F}}

\def\SS{{\cal S}}
\def\OO{{\cal O}}

\def\com{\ar@{}[rd]|{\circlearrowleft}}

\title {Exemples de Lattès et domaines faiblement sphériques de $\C^n$}
\author {Christophe Dupont}
\date{}

\begin{document}

\maketitle

\begin{abstract}
Les bassins d'attraction des relevés polynomiaux d'exemples de Lattès forment une classe remarquable de domaines pseudoconvexes bornés. Nous montrons que leur frontière s'obtient comme quotient d'une hypersurface sphérique compacte, et en décrivons précisément les singularités. Ces domaines, qui sont donc naturellement apparentés à la boule euclidienne, le sont aussi au polydisque, en tant que bassins d'attractions d'applications dont le courant de Green est ``régulier''. Les fonctions th\^eta jouent un r\^ole clef dans notre approche et nous permettent aussi d'expliciter des exemples de Lattès en dimension $2$. \\

\footnotesize{\noindent 2000 Mathematics Subject Classification : 14K25, 32S25, 32T99, 32H50  

}

\end{abstract}

\section{Introduction et résultats}

En 1918, Lattès découvrit les premiers exemples de fractions rationnelles dont l'ensemble de Julia remplit toute la sphère de Riemann. Sa construction, généralisée aux espaces projectifs de toute dimension, conduit à définir un exemple de Lattès de $\Pj^k$ comme un endomorphisme holomorphe $f$, de degré $d \geq 2$, faisant commuter un diagramme :
\[
\xymatrix{
          A^k \ar[r]^D \ar[d]_\si & A^k \ar[d]^\si \\
          \Pj^k \ar[r]^f & \Pj^k
        }
\]
où $A^k$ est un tore complexe, $D$ un endomorphisme holomorphe et $\si$ un rev\^etement galoisien de groupe $G$. Dans \cite{BL2}, ces endomorphismes sont caractérisés par la régularité et la stricte positivité de leur courant de Green sur un ouvert de $\Pj^k$.\\

L'essentiel de cet article est consacré à la description précise du bord du bassin d'attraction des relevés polynomiaux à $\C^{k+1}$ d'exemples de Lattès. Notre principal résultat est le suivant :
\begin{thm}
Soit $f$ un exemple de Lattès de $\Pj^k$ induit par un rev\^etement galoisien $\si : A^k \to \Pj^k$ de groupe $G$. Soit $T$ le courant de Green de $f$ et $\Om_F$ le bassin d'attraction d'un de ses relevés polynomiaux $F$. Alors il existe un fibré en droites $L(H,\al)$ sur $A^k$ de forme hermitienne $H$ définie négative, une hypersurface sphérique compacte $\Sigma$ dans $L(H,\al)$ et un rev\^etement galoisien $\tilde \si : L(H,\al)^- \to \C^{k+1} \setminus \lbrace 0 \rbrace$ induisant $\si$ sur les bases tels que $\tilde \si(\Sigma) = \partial \Om_F$ et $\si^* T = - { \pi \over 2 } dd^c H$.

Le bord de $\Om_F$ est sphérique au voisinage de $z_0 = \tilde \si \{ x_0,u_0 \}$ si le stabilisateur $K$ de $\dot x_0$ sous l'action de $G$ est trivial. Sinon, il a pour équation :   
\[  \lbrace (y,w) \in V_0 \times (\C,0), \Re (w) - H(\Phi ^{-1}(y),\Phi ^{-1}(y)) = 0 \rbrace  \]
où $V_0$ est un voisinage de l'origine dans $\C^k$ et $\Phi : \C^k/K \to \C^k$ est un biholomorphisme dont les coordonnées forment une base de l'algèbre des polyn\^omes invariants sous l'action du groupe des parties linéaires de $K$.
\end{thm}

D'après cette description, $\Om_F$ est un domaine disqué dont le bord est analytique réel (et m\^eme sphérique) en dehors d'un très petit lieu singulier. Ainsi les applications $F : \Om_F \to \Om_F$ sont des exemples insolites d'auto-applications holomorphes propres non injectives. En effet, ni les domaines bornés strictement pseudoconvexes (cf \cite{DF}, \cite{Pi}) (voire, comme on le conjecture, tous les domaines bornés lisses de $\C^{k+1}$), ni les domaines de Reinhardt dont une portion du bord est strictement pseudoconvexe (cf \cite{B}) n'admettent de telles applications.

Donnons un aperçu de la preuve du théorème. La désingularisation sphérique s'obtient en relevant le diagramme commutatif définissant $f$ aux fibrés en droites $\OO(-1)^- \simeq  \C^{k+1} \setminus \lbrace 0 \rbrace$ et $\si ^* \OO(-1)^- \simeq L(H,\al)^-$. Cela fait l'objet de la section \ref{desing}, où l'on construit un diagramme cubique (proposition \ref{propdesing}). Cette approche, qui nous a été suggérée par J.J. Loeb, peut reposer sur le théorème d'Appell-Humbert en vertu duquel $\si^*\OO(-1)$ est isomorphe à $L(H,\al)$. Il est cependant plus simple, techniquement et conceptuellement, d'expliciter un morphisme de $L(H,\al)$ sur $\OO(-1)$ induisant $\si$ sur les bases. On utilise pour cela le fait que les coordonnées de $\si$ sont des fonctions th\^eta de m\^eme type $(H,\al)$ (cf section \ref{sec2}). La description des singularités (cf section \ref{desc}) fait appel à la théorie des invariants : il s'agit d'interpr\^eter $\tilde \si$ comme le passage au quotient par un groupe fini d'automorphismes du fibré $L(H,\al)$. Lorsque $k=1$, Ueda avait traité un cas particulier ``à la main'' dans \cite{U1} et l'analogue de notre description avait été obtenue dans \cite{BL1}, par d'autres méthodes.\\

Dans la dernière section de l'article, on exhibe des exemples de Lattès en dimension $2$, gr\^ace à la liste des couples $(A^2,G)$ vérifiant $A^2 / G=\Pj^2$ établie dans \cite{KTY}. En particulier, les endomorphismes étudiés par Ueda (cf \cite{U1}, \cite{U2}) sont en fait des exemples de Lattès :
\begin{prop}
Les endomorphismes critiquement finis suivants sont des exemples de Lattès :
\begin{displaymath}
\begin{array}{rccc}
                  & \Pj^2  &  \longrightarrow    &  \Pj^2 \\
        f_1:       & [x:y:z]  & \longmapsto & [(-x+y+z)^2:(x-y+z)^2:(x+y-z)^2] \\       
        f_2:       & [x:y:z]  & \longmapsto & [(x-y+z)^2:(-x+y+z)^2:(x+y-z)^2] \\
        f_3:       & [x:y:z]  & \longmapsto & [(x+y-z)^2:(-x+y+z)^2:(x-y+z)^2] 
\end{array}                                      
\end{displaymath}
Ils sont semi-conjugués aux dilatations :
\[ D_1 = \sqrt2 \left(
\begin{array}{cc}
     e ^{ {i\pi \over 4}}  & 0 \\
                   0       & e^{ {i\pi \over 4}}  
\end{array}
                \right),\ 
 D_2 = \sqrt2 \left(
\begin{array}{cc}
    {1 \over \sqrt 2} &   {1 \over \sqrt 2} \\
    {1 \over \sqrt 2} &  {-1 \over \sqrt 2}
\end{array}
                \right),\ 
D_3 = \sqrt2 \left(
\begin{array}{cc}
    {1 \over \sqrt 2} &   {1 \over \sqrt 2} \\
    {i \over \sqrt 2} &  {-i \over \sqrt 2}
\end{array}
                \right)\]
par le rev\^etement galoisien $A_i \times A_i \to \Pj^2$ de groupe $\langle G(4,2,2) ,{{1+i}\over2} \binom{1}{1} \rangle$.  
\end{prop}
En revanche, l'application critiquement finie suivante, étudiée par Fornaess et Sibony (\cite{FS}, §$8$), n'est pas un exemple de Lattès, bien que son ensemble de Julia soit $\Pj^2$ : 
\begin{prop} 
L'endomorphisme critiquement fini :
\begin{displaymath}
g: 
\begin{array}{rccc}
             &  \Pj^2    &  \longrightarrow    &  \Pj^2 \\
             &  [x:y:z]  & \longmapsto       & [(x-2y)^2:(x-2z)^2:x^2]  
      
\end{array}                                      
\end{displaymath}
n'est pas un exemple de Lattès.
\end{prop}
F. Berteloot et J.J. Loeb ont manifesté de l'intér\^et pour ce travail. Je les en remercie.

\section{Fibrés en droites sur les espaces projectifs et les tores complexes}\label{sec2}

\subsection{Généralités} \label{gene}
On consultera \cite{GH} ou \cite{D} pour les généralités concernant les fibrés en droites. Cette partie est destinée à fixer les notations et rappeler des propriétés bien connues.

Soient $L,L'$ deux fibrés en droites sur $X$ et $Y$, et $f:X \to Y$ une application holomorphe. Un morphisme homogène de $L$ dans $L'$ de degré $d$, induisant $f$ sur les bases, est une application holomorphe dont l'expression au dessus d'ouverts de trivialisation $U_\al, V_{\al'}$ s'écrit : 
\begin{displaymath}
\begin{array}{ccc}
                          U_\al \times \C   & \longrightarrow  &    V_{\al'} \times \C                 \\
                             (x,t)       & \longmapsto      &  (f(x),d_{\al \al'}(x).t^d)
\end{array}
\end{displaymath}
où $d_{\al \al'}$ est une fonction holomorphe non nulle. Un morphisme de fibré est un morphisme homogène de degré $1$. Le principe du maximum entra\^ine :

\begin{lem} \label{propor}
Soit $L$ un fibré en droites sur une variété complexe compacte connexe $X$, et $u_i:L \to L, i=1,2$, deux morphismes homogènes de m\^eme degré, induisant $f$ sur $X$. Alors il existe $c \in \C^*$ tel que $u_2=c.u_1$. 
\end{lem}

Soient $f:X \to Y$ une application holomorphe, $L$ un fibré en droites sur $Y$ et $f^*L$ son image réciproque. Soit $\hat f_L:f^*L \to L$ le morphisme de fibré induisant $f$ sur la base et l'identité dans les fibres. On a $\widehat {g \circ f}_L = \hat g_L \circ \hat f _{g^* L}$ et la propriété universelle suivante :
\begin{lem}\label{univ}
Soit $L$ un fibré en droites sur $X$. Avec les notations précédentes, tout morphisme de fibré $v:L \to L'$ induisant $f$ sur les bases se factorise par $\hat f_{L'}$ : il existe un isomorphisme de fibrés $\eta : L \to f^* L'$ vérifiant  $v=\hat f_{L'} \circ \eta$. 
\end{lem}

\proof
Le morphisme $v$ s'écrit en coordonnées :
\begin{displaymath}
\begin{array}{ccc}

                U_\al \times \C                   & \longrightarrow  &    V_{\al'}  \times \C            \\
                (x,t)                             & \longmapsto      &  (f(x),d_{\al \al'}(x).t)

\end{array}
\end{displaymath}
Puisque les fonctions de transition de $f^*L'$ sont celles de $L'$ composées par $f$, on définit un morphisme de fibré $\eta : L \to f^*L'$ en posant :
\begin{displaymath}
\begin{array}{ccc}

                U_\al  \times \C          & \longrightarrow  &    f^{-1}(V_{\al'})  \times \C    \\
                (x,t)                     & \longmapsto      &    (x,d_{\al \al'}(x).t)

\end{array}
\end{displaymath}
C'est un isomorphisme vérifiant $v=\hat f_{L'} \circ \eta$.\fin

\subsection{Fibrés en droites sur $\Pj ^k$}\label{fibproj}

Soient $\pi : \C^{k+1} \setminus \lbrace 0 \rbrace \to \Pj^k$ la projection canonique et $V_j:=\lbrace \pi(z) \in \Pj^k, z_j \neq 0 \rbrace$. On regarde le fibré en droites $\OO(n)$ comme la variété $\C ^{k+1} \setminus \lbrace 0 \rbrace \times \C$ quotientée par la relation d'équivalence $\sim_n$ dont les classes sont $[z,v]_n := \lbrace (\la z, \la^n v), \la \in \C \rbrace$. L'application suivante est un biholomorphisme :
\begin{displaymath}
\Psi :
\begin{array}{rccc}
                       &  \OO(-1)^-      & \longrightarrow   &   \C^{k+1} \setminus \lbrace 0 \rbrace   \\
                       &    \lbrack z,v \rbrack_{-1}     & \longmapsto      &      v \cdot z    
\end{array}
\end{displaymath}
Son inverse est donné par $\Psi^{-1}(z) =  \lbrack z,1 \rbrack_{-1}$.     

\begin{lem} \label{struc1}
Soit $f$ un endomorphisme holomorphe de $\Pj^k$ de degré $d$, et $F$ un de ses relevés polynomiaux à $\C^{k+1}$. Alors

\begin{enumerate}

\item
$f^* \OO (-1)$ est isomorphe à $ \OO (-1) ^{\otimes d} $, c'est à dire à $\OO (-d)$. 

\item
Les morphismes homogènes $\OO (-1) \to \OO (-1) $ de degré $d$ induisant $f$ sur $\Pj ^k$ sont de la forme $[z,v]_{-1} \mapsto [F(z),c.v^d]_{-1}$, où $c \in \C^*$. 

\end{enumerate}
\end{lem}

\proof
Commençons par le premier point. Gr\^ace à la propriété universelle du fibré image réciproque (cf lemme \ref{univ}), il suffit d'exhiber un morphisme de fibré $\OO (-d) \to  \OO (-1)$ induisant $f$ sur $\Pj^k$. Puisque $F$ est homogène de degré $d$, l'application :
\begin{displaymath}
\begin{array}{ccc}

  \C ^ {k+1} \setminus \lbrace 0 \rbrace \times \C  & \longrightarrow  &   \C ^ {k+1} \setminus \lbrace 0 \rbrace \times \C  \\
                (z,v)                               & \longmapsto      &  (F(z),v)

\end{array}
\end{displaymath}
passe au quotient pour $\sim_{-d}$ à la source et $\sim_{-1}$ au but, ce qui permet de conclure. Pour le second point, on utilise l'homogénéité de $F$ pour vérifier que l'application : 
\begin{displaymath}
\begin{array}{ccc}

  \C ^ {k+1} \setminus \lbrace 0 \rbrace \times \C  & \longrightarrow  &   \C ^ {k+1} \setminus \lbrace 0 \rbrace \times \C  \\
                (z,v)                               & \longmapsto      &  (F(z),v^d)

\end{array}
\end{displaymath}
passe au quotient pour $\sim_{-1}$ à la source et au but. On utilise alors le lemme \ref{propor}.
\fin

\subsection{Fibrés en droites sur un tore complexe $A^k$} \label{fibtore}

Soit $\Ga$ un réseau de $\C^k$ tel que $A^k = \C^k / \Ga$. On notera $\Pi : \C^k \to  \C^k / \Ga$ la projection canonique et $\dot x$ pour $\Pi(x)$. Un type est un couple $(H,\al)$ où

\noindent - $H$ est une forme hermitienne sur $\C^k$, linéaire à droite, telle que $\Im H (\Ga,\Ga) \subset \Z$

\noindent - $\al$ est une application de $\Ga$ dans $\Sg^1$, vérifiant $\al(\ga_1+\ga_2)=\al(\ga_1)\al(\ga_2)(-1)^{\Im H (\ga_1,\ga_2)}$

On définit une action de $\Ga$ sur $\C^k \times \C$ par $\ga \cdot (x,u)  := (x+\ga,e_\ga (x).u)$, où $e_\ga (x) := \al(\ga).e ^{ \pi [H (\ga,x)+ {1 \over 2} H (\ga,\ga)]}$. On note $L(H,\al)$ le fibré en droites sur $A^k$ obtenu comme le quotient de $\C^k \times \C$ par cette action, $\lbrace x,u \rbrace _{(H,\al)}$ ses éléments et $q$ la métrique induite par la fonction  $q(x,u)=e ^ {- { \pi \over 2} H(x,x)} \abs u$. Le produit $L(H_1,\al_1) \otimes  L(H_2,\al_2)$ est isomorphe à $L(H_1+H_2,\al_1.\al_2)$. L'espace des sections de $L(H,\al)$ est isomorphe à l'ensemble des fonctions th\^eta normalisées de type $(H,\al)$, auquel on ajoute la fonction nulle. On rappelle que ce sont les fonctions holomorphes sur $\C^k$, non identiquement nulles, vérifiant $\theta(z+\ga)=e_\ga(x).\theta(z)$. Puisqu'une telle fonction est de type $(0,1)$ dès qu'elle ne s'annule pas, on a $L(H_1,\al_1) \simeq L(H_2,\al_2)$ si et seulement si $(H_1,\al_1) = (H_2,\al_2)$. Le fait classique suivant (cf \cite{D}, Chap IV, §3) jouera un r\^ole essentiel dans la démonstration : 

\begin{thm}\label{the}
Toute application holomorphe $\si : A^k \to \Pj^k$ est induite par une application $\theta = (\theta_0,...,\theta_k)$ dont les coordonnées sont des fonctions th\^eta normalisées de m\^eme type. 
\end{thm}
Une preuve identique à celle du lemme \ref{struc1} fournit : 
\begin{lem}
Soit $\si : A^k \to \Pj^k$ induite par $\theta = (\theta_0,...,\theta_k)$, où les fonctions th\^eta sont de type $(H,\al)$. Les morphismes de fibrés $\phi : L(dH,\al^d) \to \OO(d)$ induisant $\si$ sur les bases sont de la forme $\phi \{x,u\}_{(dH,\al^d)} = [ \theta (x), c.u ]_d$, où $c \in \C^*$. 
\end{lem}
Notons qu'alors $\si ^* \OO(d) \simeq L(H,\al)$ par le lemme \ref{univ}. On donne maintenant l'analogue du lemme \ref{struc1}. Soit $\vphi = \vec \vphi + \tau$ un endomorphisme de $A^k$. On note $\tilde \tau$ l'application de translation de vecteur $\tau$ sur $A^k$, de sorte que $\vphi=\tilde \tau \circ \vec \vphi$ et $\vphi^*L(H,\al)=\vec \vphi^* \, \tilde \tau ^* \, L(H,\al)$. Etant donné un type $(H,\al)$, on définit $H_\vphi$ et $\al_\vphi$ par :
\begin{displaymath}
\begin{array}{rcl}
         H_\vphi (w,w') &:= &H(\vec \vphi w,\vec \vphi w' ) \\
         \al_\vphi (\ga)&:= &\al(\vec \vphi \ga) e ^{ 2i\pi \Im H(\vec \vphi \ga,\tau)}
\end{array}
\end{displaymath} 
Notons que $\al_\vphi$ est bien définie, puisque $\vec \vphi \left( \Ga \right) \subset \Ga$, et prend ses valeurs dans $\Sg^1$. L'identité $\al_\vphi (\ga_1+\ga_2)= \al_\vphi (\ga_1) \al_\vphi (\ga_2) (-1)^{\Im H_\vphi  (\ga_1,\ga_2)}$ montre que $(H_\vphi ,\al_\vphi )$ est encore un type, dont on note $(e_\ga^\vphi )_{\ga \in \Ga}$ les multiplicateurs.

\begin{lem} \label{struc2}
Soit $\vphi = \vec \vphi + \tau$ un endomorphisme de $A^k$. Alors :

\begin{enumerate}

\item
$\vphi ^*L(H,\al)$ est isomorphe à $L(H_\vphi ,\al_\vphi )$

\item
Si $(H_\vphi ,\al_\vphi )=(dH,\al^d)$, on a $e ^{ \pi H (\tau , \vec \vphi \ga) }.e_\ga^d (x) = e_{\vec \vphi \ga} (\vphi x)$.

\item
Il existe des morphismes homogènes $L(H,\al) \to L(H,\al)$ de degré $d$ induisant $\vphi $ sur $A ^k$ si et seulement si $(H_\vphi ,\al_\vphi )=(dH,\al^d)$. Ils sont alors de la forme $\lbrace x,u \rbrace_{(H,\al)} \mapsto  \lbrace \vphi x, c.e ^{ \pi H (\tau , \vec \vphi x) }.u^d \rbrace_{(H,\al)}$, où $c \in \C^*$.

\end{enumerate}
\end{lem}

\proof
Pour le premier point, on vérifie que l'application : 
\begin{displaymath}
\chi_1 :
\begin{array}{cccc}
        &\C ^ k \times \C   & \longrightarrow  &    \C ^ k \times \C  \\
        & (x,u)            & \longmapsto      &  (\vphi x , e ^{ \pi H(\tau,\vec \vphi x)} . u)
\end{array}
\end{displaymath}
induit un morphisme de fibré entre $L(H_\vphi ,\al_\vphi)$ et $L(H,\al)$, puis on applique le lemme \ref{univ}. On a : \[ \chi_1 \left( x+\ga, e_\ga ^\vphi (x).u \right) =  (\vphi x + \vec \vphi \ga ,  e_\ga ^\vphi (x).e ^{ \pi H(\tau,\vec \vphi \ga)}.e ^{ \pi H(\tau,\vec \vphi x)}.u). \]
Il suffit donc de vérifier que $e_\ga ^\vphi (x).e ^{ \pi H(\tau,\vec \vphi \ga)} =  e_{\vec \vphi \ga} (\vphi x)$ : 
\begin{displaymath}
\begin{array}{rcl}
   e_\ga ^\vphi (x).e ^{ \pi H(\tau,\vec \vphi \ga)}  &  =  &   \al(\vec \vphi \ga).e ^{ 2i\pi \Im H (\vec \vphi \ga,\tau)}.e ^{ \pi [H (\vec \vphi \ga,\vec \vphi x) + {1 \over 2} H (\vec \vphi \ga, \vec \vphi \ga)]}.e ^{ \pi H(\tau,\vec \vphi \ga)} \\
                                                      &  =  &   \al(\vec \vphi \ga).e ^{ \pi [H (\vec \vphi \ga,\vec \vphi x) + {1 \over 2} H (\vec \vphi \ga, \vec \vphi \ga)]}.e ^{ \pi H(\vec \vphi \ga, \tau)} \\
                                                      &  =  &    e_{\vec \vphi \ga} (\vphi x)
\end{array}
\end{displaymath}
Pour le deuxième point, on vérifie que l'hypothèse entra\^ine immédiatement :
\begin{displaymath}
\begin{array}{rcl}
  e ^{ \pi H (\tau , \vec \vphi \ga) }.e_\ga^d(x)  &  =  &    e ^{ \pi H (\tau , \vec \vphi \ga)}.\al(\vec \vphi \ga).e ^{ 2i\pi \Im H (\vec \vphi \ga,\tau)}.e ^{ \pi [H (\vec \vphi \ga,\vec \vphi x) + {1 \over 2} H (\vec \vphi \ga, \vec \vphi \ga)]}     \\                                                      
                                                        &  =  &   \al(\vec \vphi \ga).e ^{ \pi [H (\vec \vphi \ga,\vec \vphi x + \tau) + {1 \over 2} H (\vec \vphi \ga, \vec \vphi \ga)]} = e_{\vec \vphi \ga} (\vphi x)
\end{array}
\end{displaymath}
Montrons le dernier point. Lorsque $d=1$, la propriété universelle du fibré image réciproque (cf lemme \ref{univ}), l'unicité du type et le premier point entra\^inent la condition nécessaire. Lorsque $d \geq 2$, on remarque que tout morphisme homogène $L \to L$ de degré $d$ donne naissance à un morphisme de fibré $L^{\otimes d} \to L$ lorsque l'on supprime la puissance $d$ dans les fibres. Pour la condition suffisante et les expressions annoncées, on vérifie que l'application :
\begin{displaymath}
\chi_2 :
\begin{array}{cccc}
    &\C ^ k \times \C   & \longrightarrow  &    \C ^ k \times \C  \\
    & (x,u)            & \longmapsto      &   \left( \vphi x, e ^{ \pi H (\tau , \vec \vphi x) }.u^d \right)
\end{array}
\end{displaymath}
passe au quotient pour $\sim_{(H,\al)}$ à la source et au but. Le deuxième point permet de conclure puisqu'on a :
\[ \chi_2 \left( x+\ga, e_\ga(x).u \right) = \left( \vphi x + \vec \vphi \ga , e ^{ \pi H (\tau , \vec \vphi x + \vec \vphi \ga) }.e_\ga^d (x). u^d \right) \]
Cela achève la démonstration du lemme.\fin

\ 
 
On termine cette partie en explicitant des trivialisations du fibré $L(H,\al)$. Fixons $U$ un ouvert de $\C^k$, disjoint de ses translatés par l'action de $\Ga$ ($U$ sera dit $\Ga$-petit). L'application $\Ga_U : \bigcup_{\ga \in \Ga} \left( U+\ga \right) \to \Ga$ qui à un point $x$ associe l'unique élément $\ga \in \Ga$ vérifiant $x+\ga \in U$ est alors bien définie. Soit $\epsilon$ une fonction holomorphe ne s'annulant pas sur $U$. L'application suivante :
\begin{displaymath}
\begin{array}{ccc}

  \bigcup_{\ga \in \Ga} \left( U+\ga \right) \times \C       & \longrightarrow  &   U \times \C    \\
          (x,u)                                              & \longmapsto      &  \left( x+\Ga_U (x),\epsilon(x+\Ga_U (x)). e_{\Ga_U (x)}(x).u \right) 

\end{array}
\end{displaymath}
passe au quotient modulo $\sim_{(H,\al)}$ et induit une trivialisation  $\psi_{(U,\epsilon)}$ de $L(H,\al)$ : 
\begin{displaymath}
\psi_{(U,\epsilon)}:
\begin{array}{rccccc}
       &  p^{-1} \circ \Pi (U)                & \longrightarrow     & \Pi(U) \times \C \\
       & \lbrace x,u \rbrace _{(H,\al)}       & \longmapsto         & \left( \Pi(x),\epsilon(x+\Ga_U (x)). e_{\Ga_U (x)}(x).u \right)
\end{array}
\end{displaymath}

\section{Désingularisation sphérique du bord du bassin d'attraction}\label{desing}

Soit $f: \Pj^k  \to \Pj^k$ un exemple de Lattès de degré $d \geq 2$ et $F: \C^{k+1} \to \C^{k+1}$ un relevé polynomial de $f$. On note $G_F$ la fonction $\lim_{p \to +\infty} { 1 \over d^p } . \log \bignorm { F^p }$, $T$ le courant de Green (donné par $\pi ^* T := dd^c G_F$) et $\Om_F$ le bassin d'attraction de $F$. On consultera \cite{HP} et \cite{S} pour plus de précisions. Par définition, on dispose d'un diagramme : 
\[
\xymatrix{
          A^k \ar[r]^D \ar[d]_\si & A^k \ar[d]^\si \\
          \Pj^k \ar[r]^f & \Pj^k
        }
\]
où $\si$ est une application de passage au quotient par un groupe $G$. On confondra $D$ et un de ses relevés $\vec D + \tau$ à $\C^k$. Soit  $(-H,\al^{-1})$ le type des fonctions th\^eta normalisées qui induisent $\si$. On a alors $\si^* \OO(-1) \simeq L(H,\al)$. On note $\tilde \theta$ le morphisme de fibré :
\begin{displaymath}
 \tilde \theta :
\begin{array}{cccc}
                 & L(H,\al)                & \longrightarrow  &    \OO(-1)     \\
                 & \{ x,u \}_{(H,\al)}     & \longmapsto      &  [ \theta (x), u ]_{-1}
\end{array}
\end{displaymath}
La forme $-H$ est définie positive puisque $\si$ est à fibres finies (cf \cite{D}, Chap.IV, Cor.$3.5$). Ainsi, une surface de niveau de $q$ est localement biholomorphe à un ouvert de la sphère unité de $\C^{k+1}$ (on dira qu'elle est sphérique). En effet : soit $\lbrace x_0,u_0 \rbrace_{(H,\al)} \in \lbrace q=c \rbrace$, $U_{x_0}$ un voisinage de $x_0$ dans $\C^k$ et $\log$ une détermination du logarithme en ${u_0 \over c} \neq 0$. Si $v = {2 \over \pi}  \log ({u \over c})$, l'équation de $\lbrace q=c \rbrace$ s'écrit :
\[ \lbrace (x,v) \in U_{x_0} \times (\C,v_0) \, , \,  \Re (v) - H(x,x) = 0 \rbrace \]
On reconna\^it l'équation d'un ouvert de la sphère unité de $\C^{k+1}$, dans sa version non bornée (cf \cite{R}, Chap.$2$, §$3$). Une telle hypersurface est compacte puisqu'elle est fibrée en cercles au dessus de $A^k$. Le lemme qui suit est crucial.  
\begin{lem} \label{herm}
$ $

\begin{enumerate}

\item
$(H_D,\al_D)=(dH,\al^d)$

\item
Pour tout $g \in G$, $(H_g,\al_g) = (H,\al)$. 
\end{enumerate}
\end{lem}

\proof
Montrons le premier point. La relation $\si \circ D = f \circ \si$ entra\^ine $D^* \si^* \OO(-1) = \si^* f^* \OO(-1)$. Le lemme \ref{struc1} implique :
\[ {\si^* f^* \OO(-1) \simeq \si^* \OO(-1)^{\otimes d}  \simeq  L(H,\al)^{\otimes d}  \simeq  L(dH,\al^d) } \]
D'autre part, le lemme \ref{struc2} donne :
\[ { D^* \si^* \OO(-1) \simeq D^* L(H,\al) \simeq L(H_D,\al_D) } \]
L'unicité du type permet de conclure. On montre de m\^eme le second point.\fin

\ 

La proposition suivante désingularise $\partial \Om_F$ en une hypersurface sphérique :
\begin{prop}\label{propdesing}
Il existe $\DD : L(H,\al) \to L(H,\al)$ un morphisme homogène de degré $d$ et $F$ un relevé polynomial de $f$ tel que le diagramme suivant commute :   
\[
\xymatrix{
          L(H,\al)^- \ar[rr]^{\DD} \ar[dd]_{\tilde\si} \ar[dr] && L(H,\al)^- \ar[dr] \ar[dd]^{\tilde\si} \\
          & A^k \ar[rr]^{\hskip -1 cm D} \ar[dd]_\si && A^k \ar[dd]^\si \\
          \C^{k+1} \setminus \lbrace 0 \rbrace  \ar[rr]^{\hskip 1 cm F}  \ar[dr]_\pi && \C^{k+1} \setminus \lbrace 0 \rbrace  \ar[rd]_\pi \\
          & \Pj^k \ar[rr]^f && \Pj^k \\
        }
\]
où $\tilde \si$ désigne l'application $\Psi \circ \tilde \theta$. Quitte à normaliser la métrique $q$, on a $e^ {G_F \circ \tilde \si} = q$, $\tilde \si \lbrace q = 1 \rbrace = \partial \Om_F$ et l'égalité $\si ^* T = -{\pi \over 2} dd^c H$.   
\end{prop}

\proof
Les lemmes \ref{struc1}, \ref{struc2} et \ref{herm} assurent l'existence des morphismes homogènes de degré $d$ suivant :
\begin{displaymath}
\begin{array}{rccl}
  \DD : &  \lbrace x,u \rbrace_{(H,\al)}  & \longmapsto & \lbrace D x, e ^{ \pi H (\tau , \vec D x)  }.u^d \rbrace_{(H,\al)}  \\

  \F      : &  [z,v]_{-1}               & \longmapsto &   [F(z),c.v^d]_{-1}
\end{array}
\end{displaymath}
Nous allons montrer qu'il existe $c \in \C^*$ tel que $\tilde \theta \circ \DD  = \F \circ \tilde \theta$. Cela revient à établir 
\begin{equation}\label{commutatif}
\exists c \in \C^*, \ \forall x \in \C^k, \  c.F \circ \theta (x) =  \theta \circ D(x)  . e ^{ \pi H (\tau , \vec D x) }
\end{equation}
L'identité $ \pi \circ \theta \circ D = \pi \circ F \circ \theta $ détermine une unique fonction holomorphe non nulle $c=c(x)$ vérifiant (\ref{commutatif}) sur $\C^k$. Le calcul suivant, où $e_\ga$ est de type $(H,\al)$, montre que $c$ est $\Ga$-périodique.
\begin{displaymath}
\begin{array}{rcl}
 c(x+\ga).F \circ \theta (x+\ga) & = & \theta \circ D (x+\ga). e ^{ \pi H (\tau , \vec D(x+\ga)) } \\
                                 & = &  e_{\vec D \ga}^{-1} (Dx).\theta \circ D (x) .  e ^{ \pi H (\tau , \vec D x) } . e ^{ \pi H (\tau , \vec D \ga) } \\                              
                                 & = &  e_{\vec D \ga}^{-1} (Dx).c(x).F \circ \theta (x). e ^{ \pi H (\tau , \vec D \ga) }\\
                                 & = &  e_{\vec D \ga}^{-1} (Dx).c(x).F \circ \theta (x+\ga).e_\ga^d (x). e ^{ \pi H (\tau , \vec D \ga) }\\
                                 & = & c(x).F \circ \theta (x+\ga)

\end{array}
\end{displaymath}
La dernière égalité provient du lemme \ref{struc2}. Puisque l'une des composantes de $F$ est non nulle, on a bien $c(x+\ga) = c(x)$. Cette fonction est donc constante par le théorème de Liouville. On supposera qu'elle vaut $1$ quitte à changer $F$ en $c.F$. Le diagramme annoncé provient alors de $\Psi \circ \F  = F \circ \Psi$, où $\Psi$ est le biholomorphisme entre $\OO(-1)^-$ et $\C^{k+1} \setminus \{ 0 \}$ introduit dans la partie \ref{fibproj}. 

A présent, on normalise la métrique $q$ sur $L(H,\al)$ afin d'obtenir la relation $q \circ \DD = q ^d$. On notera $q_\delta$ pour $\delta.q$. L'égalité $H_D=dH$ (cf lemme \ref{herm}) entra\^ine :
\begin{displaymath}
\begin{array}{rcl}

 q_\delta \circ \DD \lbrace x,u \rbrace & = &   \delta . e ^ {- { \pi \over 2} H(\vec Dx + \tau,\vec Dx + \tau)} \abs {e ^{ \pi H (\tau , \vec D x)} . u^d}  \\

                                           & = &   \delta . e ^ {- { \pi \over 2} [H(\vec Dx,\vec Dx) + 2 \Re H(\vec D x,\tau) + H(\tau,\tau)]} . e ^{ \pi \Re H (\tau , \vec D x)} . \abs {u^d}  \\    

                                           & = &  \delta . e ^ {- { \pi \over 2} dH(x,x)} . e ^ {- { \pi \over 2} H(\tau,\tau)} . \abs {u^d}
    
\end{array}                                      
\end{displaymath}
Cette expression coïncide avec $q _\delta ^d \lbrace x,u \rbrace = \left( \delta.e ^ {- { \pi \over 2} H(x,x)}.\abs u \right) ^d $ si $\delta$ vérifie $\delta.e ^ {- { \pi \over 2} H(\tau,\tau)} = \delta ^d$, c'est à dire si $\delta = \left( e ^ {- { \pi \over 2} H(\tau,\tau)} \right)^ {1 \over {d-1}}$.

Montrons l'égalité $e^ {G_F \circ \tilde \si} = q$. Soit $\bignorm{.}$ une norme sur $\C^{k+1}$ : puisque $A^k$ est une variété compacte, il existe une constante $C$ telle que ${1 \over C} .q \leqslant \bignorm{ \tilde \si } \leqslant C.q$. D'après ce qui précède, on a $q \circ \DD^p  =  q ^{d^p}$ et $\tilde \si \circ \DD^p  =  F^p \circ \tilde \si$ pour tout $p$. Ainsi, $G_F \circ \tilde \si = \lim_{p \to +\infty} { 1 \over d^p } . \log \bignorm { F^p \circ \tilde \si } = \lim_{p \to +\infty} { 1 \over d^p } . \log \bignorm { \tilde \si \circ \DD^p } = \lim_{p \to +\infty} { 1 \over d^p } . \log q \circ \DD^{d^p} = \log q$. L'expression de $\tilde \si$ et l'homogénéité de $G_F$ entra\^inent $G_F \circ \theta (x) =  \log \delta - { \pi \over 2} H(x,x)$ et $\si ^* T = -{\pi \over 2} dd^c H$ se déduit de l'identité $\pi \circ \theta = \si \circ \Pi$.\fin

\ 

\begin{rem}
\rm L'égalité $\si ^* T = -{\pi \over 2} dd^c H$ permet de retrouver le caractère défini positif de $-H$. En effet, $-H$ est positive puisque $T$ est un courant positif. Si cette forme était dégénérée, $T$ aurait un potentiel local maximal en dehors de $\si(Crit \; \si)$. La mesure de probabilité $\mu=T^k$ serait alors portée par $\si(Crit \; \si)$, ce qui est impossible d'après l'inégalité de Chern-Levine-Nirenberg (cf \cite{S}, Prop.A.$6.3$).
\end{rem}

\section{Description du bord du bassin d'attraction}\label{desc}

Notons $\CC := \{   \lbrace x,u \rbrace_{(H,\al)}, \dot x \in Crit(\si) \} $ l'ensemble critique de $\tilde \si$, et $\SS$ ses valeurs critiques, de sorte que $\tilde \si : L(H,\al) \setminus \CC \to \C^{k+1} \setminus \SS$ soit un rev\^etement fini. Puisque le bord de $\Om_F$ est égal à $\tilde \si \lbrace q=1 \rbrace$, il est sphérique en dehors de $\SS$. Sa description aux points de $\SS$ nécessite l'interprétation de $\tilde \theta$ suivante :

\begin{lem}
L'application $\tilde \theta$ s'identifie au passage au quotient par un groupe fini $\tilde G$ d'automorphismes de $L(H,\al)$ et son action laisse invariants les niveaux de $q$.
\end{lem}

\proof
D'après le lemme \ref{univ}, il existe un isomorphisme $\eta :  L(H,\al) \to \si^* \OO(-1)$ vérifiant $\tilde \theta = \hat \si _{\OO(-1)} \circ \eta$. Par définition, $\hat \si_{\OO(-1)}$ est l'application de passage au quotient de $\si^*\OO(-1)$ par le groupe $\hat G  := \{ \hat g_{\si^*\OO(-1)} : \si ^* \OO(-1) \to \si^* \OO(-1) , g \in G \}$. Ainsi, $\tilde \theta$ est celle de $L(H,\al)$ par le groupe conjugué $\tilde G := \eta ^{-1} \circ \hat G \circ \eta$. L'égalité $e^ {G_F \circ \Psi \circ \tilde \theta} = q$ de la proposition \ref{propdesing} montre que les surfaces de niveau de $q$ sont invariantes par $\tilde G$.
\fin

\ 

Le bord de $\Om_F$ s'identifie alors au quotient de $\lbrace q=1 \rbrace$ par $\tilde G$ ($\Psi$ est un biholomorphisme). Fixons $\lbrace x_0,u_0 \rbrace_{(H,\al)} \in \CC$. On note $G_0$ le stabilisateur de $\dot{x_0}$ et $\vec G_0$, $\tilde G_0$, les sous groupes de $\vec G$, $\tilde G$ correspondant. On suppose que $\vec G$ est un sous groupe de $\U_k(\C)$ quitte à effectuer un changement linéaire de coordonnées. Soit $U_{x_0}=B(x_0,\nu)$ une boule euclidienne $\Ga$-petite, telle que les translatées $g\Pi(U_{x_0})$ soient disjointes lorsque $g \notin G_0$. Vérifions que $\Pi(U_{x_0})$ est invariante par $G_0$ (on reprend les notations de la partie \ref{fibtore}). Prenons $g := \vec g + \kappa \in G_0$, et $\ga := \Ga_{U_{x_0}}(g x_0)$, de sorte que $g(x_0+t)=\vec g x_0 + \vec g t + \kappa = x_0 + \vec g t - \ga$, pour $t \in B(0,\nu)$. On a alors $g(U_{x_0}) = U_{x_0} - \ga$ ($\vec g \in \U_k(\C)$). Observons aussi que $\Ga_{U_{x_0}}(g (x_0+t)) = \Ga_{U_{x_0}}(g x_0) = \ga$. 

Etudier le bord de $\Om_F$ en $z_0 := \tilde \si \lbrace x_0,u_0 \rbrace_{(H,\al)}$ revient alors à étudier le quotient de $\lbrace q=1 \rbrace$ au dessus de $\Pi(U_{x_0})$ par $\tilde G_0$. Pour cela, on exhibe des trivialisations dans lesquelles $\tilde G_0$ induit l'identité dans les fibres. On dit qu'elles sont $\tilde G_0$-équivariantes. 
\begin{lem}
La trivialisation $\psi_{(U_{x_0},\epsilon)}$ est $\tilde G_0$-équivariante si et seulement si $\epsilon$ est de la forme $\epsilon_0 \cdot m$, où $\epsilon_0$ est la fonction :
\begin{displaymath}
 \epsilon_0 :
\begin{array}{rccc}
                &             U_{x_0}      & \longrightarrow &  \C  \\ 
                &             x_0+t        & \longmapsto     &  e ^ {- \pi H(x_0,t) } 
\end{array}                                      
\end{displaymath}  
et $m: t \mapsto m(x_0+t)$ est holomorphe non nulle sur $B(0,\nu)$ et $\vec G_0$-invariante.
\end{lem}

\proof
Rappelons que $\psi_{(U_{x_0},\epsilon)}$ a pour expression (cf § \ref{fibtore}): 
\begin{displaymath}
\psi_{(U_{x_0},\epsilon)}:
\begin{array}{rccccc}
       &  p^{-1} \circ \Pi (U_{x_0})                & \longrightarrow     & \Pi(U_{x_0}) \times \C \\
       & \lbrace x,u \rbrace_{(H,\al)}              & \longmapsto         & \left( \dot x, \epsilon(x +\Ga_{U_{x_0}} (x)). e_{\Ga_{U_{x_0}} (x)}(x).u \right)
\end{array}
\end{displaymath}
Afin d'alléger l'écriture, on remplacera $\psi_{(U_{x_0},\epsilon)}$ par $\psi$. Cette trivialisation est équivariante si et seulement si on a les identités suivantes sur $p^{-1} \circ \Pi (U_{x_0})$ : 
\begin{equation}\label{condequi1}
 \forall g \in G_0, \  \psi \circ \tilde g \lbrace x_0+t,u \rbrace_{(H,\al)} =   (g(\dot{x_0} + \dot t) , \epsilon(x_0+t).u)
\end{equation}
Le lemme \ref{struc2} donne l'expression de l'automorphisme $\tilde g$ du fibré $L(H,\al)$, induit par $g = \vec g + \kappa$ sur la base. Il existe $\rho (g) \in \C$ tel que :
\begin{displaymath}
   \tilde g : 
\begin{array}{cccc}
                & L(H,\al)                                      & \longrightarrow   &                   L(H,\al)  \\
                &     \left\lbrace x,u \right\rbrace_{(H,\al)}             & \longmapsto      &    \left\lbrace g(x),\rho (g).e ^{ \pi H(\kappa,\vec g x)}.u \right\rbrace_{(H,\al)}
\end{array}
\end{displaymath}
Ainsi, l'équation (\ref{condequi1}) que l'on doit vérifier devient :
\begin{equation}\label{condequi2}
 \psi \circ \left\lbrace g(x_0+t),\rho (g).e ^{ \pi H(\kappa,\vec g (x_0+t))}.u \right\rbrace_{(H,\al)} =  (g(\dot{x_0} + \dot t), \epsilon(x_0+t).u)
\end{equation}
Comme nous l'avons observé, $\ga = \Ga_{U_{x_0}}(g (x_0+t)) = \Ga_{U_{x_0}}(g (x_0))$. Donc, par définition de $\psi$, le membre de gauche de (\ref{condequi2}) est égal à :
\[  \left(     g(\dot{x_0} + \dot t)   , e_\ga(x_0+\vec g t-\ga) . \epsilon( x_0+ \vec g t)    . \rho (g) . e ^{ \pi H(\kappa,\vec g (x_0+t))}.u \right) \]
La condition nécessaire et suffisante d'équivariance (\ref{condequi2}) s'écrit alors : 
\begin{equation}\label{condequi3}
\forall t \in B(0,\nu)\ , \ \forall g \in G_0 \ , \  {\epsilon (x_0+t) \over \epsilon( x_0+ \vec g t)}  = e_\ga (x_0+\vec g t-\ga) . \rho (g) . e ^{ \pi H(\kappa,\vec g (x_0+t))} 
\end{equation}
où $\ga =  \Ga_{U_{x_0}}(g(x_0+t)) = \Ga_{U_{x_0}}(g(x_0))$. Remarquons aussi que : 
\begin{equation}\label{sos}
 \forall g \in G_0 \ , \  e_\ga (x_0 - \ga) . \rho (g) . e ^{ \pi H(\kappa,\vec g x_0)} = 1.
\end{equation} 
En effet, puisque $\tilde \theta = \tilde \theta \circ \tilde g$, on a : (avec $e_\ga$ de type $(H,\al)$)
\begin{displaymath}
\begin{array}{rcl}
 [ \theta (x_0+t), u ]_{-1}  & = &   [ \theta ( g(x_0+t)), \rho (g) . e ^{ \pi H(\kappa,\vec g (x_0+t))}.u ]_{-1} \\
                         & = &   [ \theta ( x_0+ \vec g t - \ga), \rho (g) . e ^{ \pi H(\kappa,\vec g (x_0+t))}.u ]_{-1} \\ 
                         & = &   [ \theta ( x_0+ \vec g t ), e_{-\ga}^{-1} (x_0 + \vec g t) .\rho (g) . e ^{ \pi H(\kappa,\vec g (x_0+t))}.u ]_{-1} \\ 
\end{array}                                      
\end{displaymath}
Comme $e_{-\ga} (x_0).e_{\ga} (x_0 - \ga)=1$, on évalue en $t=0$ pour obtenir (\ref{sos}). Ainsi, le membre de droite de (\ref{condequi3}) se simplifie. Il est égal à :
\begin{displaymath}
\begin{array}{rcl}
  { \displaystyle  e_\ga (x_0+\vec g t-\ga) . \rho (g) . e ^{ \pi H(\kappa,\vec g (x_0+t))}  \over  \displaystyle e_\ga (x_0 - \ga) . \rho (g) . e ^{ \pi H(\kappa,\vec g x_0)} } & = &  { \displaystyle e ^{ \pi H (\ga,x_0 + \vec g t - \ga) } \over \displaystyle e ^{ \pi H (\ga,x_0 - \ga) }} .e ^{ \pi H (\kappa,\vec g t)}   \\        & = &   e ^{ \pi H (\ga + \kappa , \vec g t) }  \\
            & = &   e ^{ \pi H (x_0 - \vec g x_0 , \vec g t) } \\
            & = &   e ^{ - \pi H (\vec g x_0 , \vec g t) } . e ^{  \pi H (x_0 , \vec g t) }  \\
            & = &   e ^{ - \pi H (x_0 , t) } . e ^{ \pi H (x_0 , \vec g t)  }
\end{array}                                      
\end{displaymath} 
La dernière égalité provient de l'invariance de $H$ par $\vec G$ (cf lemme \ref{herm}). Donc, $\psi$ est équivariante si et seulement si : 
\begin{equation}
\forall t \in B(0,\nu)\ , \ \forall g \in G_0 \ , \  {\epsilon (x_0+t) \over \epsilon( x_0+ \vec g t)}  =  { e ^{ - \pi H (x_0 , t) } \over  e ^{ -\pi H (x_0 , \vec g t) } } 
\end{equation}  
Pour conclure, il suffit de remarquer que la fonction $\epsilon_0 (x_0+t) =  e ^{ - \pi H (x_0 , t) }$ vérifie cette équation et que le rapport de deux solutions est une fonction holomorphe non nulle $m : t \to m(x_0 +t)$ sur $B(0,\nu)$, invariante par $\vec G_0$. \fin

\ 

Désormais, nous travaillons avec la trivialisation équivariante $\psi$ donnée par $\epsilon_0$. Quitte à identifier $\Pi(U_{x_0})$ avec $U_{x_0}$ et effectuer une translation, $\psi$ prend ses valeurs dans $B(0,\nu) \times \C$. Dans ces coordonnées,  $\lbrace q=1 \rbrace$ s'écrit :
\[ \{ (t , v) \in B(0,\nu) \times (\C,v_0) \, , \,  \delta_0 .  e ^ {- { \pi \over 2} H(t,t)} . \abs v = 1 \} \]
où $\delta_0 := \delta.e ^{-{\pi \over 2} H(x_0,x_0)}$. Cela provient du calcul :
\[ q \left( x_0 + t , \epsilon_0(x_0+t)^{-1} .v  \right)   =  \delta. e ^{-{\pi \over 2} H( x_0 + t , x_0 + t )}. e ^ { \pi  \Re H( x_0,t )} \abs v  = \delta. e ^{-{\pi \over 2} H(x_0,x_0)}. e ^{-{\pi \over 2} H(t,t)}  \abs v \]
L'action de $\tilde G_0$ s'identifie à celle de $\vec G_0$ sur $B(0,\nu)$. Donc $\partial \Om_F$ en $z_0$ s'écrit  :
\[ \partial \Om_F = \lbrace (\tilde t,v) \in B(0,\nu)/{\vec G_0} \times (\C,v_0) \, , \, \delta_0 .  e ^ {- { \pi \over 2} H(\tilde t,\tilde t)} . \abs v = 1 \rbrace \]
où $t \mapsto \tilde t$ désigne la projection sur $(\C^k,0)/{\vec G_0}$. Notons que la singularité $B(0,\nu)/{\vec G_0}$, isomorphe à $\Pi(U_{x_0})/G_0$, est lisse, puisque le quotient $A^k/G$ est égal à $\Pj^k$. Le groupe $\vec G_0$ doit donc \^etre engendré par des réflexions (cf \cite{Pr}, §4), ie des éléments laissant un hyperplan fixe. L'algèbre des polyn\^omes invariants par $\vec G_0$ est alors engendrée par $k$ polyn\^omes homogènes $(P_1,...,P_k)$ algébriquement indépendants (cf \cite{ST}, Th.$5.1$), et l'application suivante est un biholomorphisme (cf \cite{C}, Th.$3$) :
\begin{displaymath}
\Phi :
\begin{array}{rccc}
           &  B(0,\nu)/{\vec G_0}            &  \longrightarrow    &    V_0           \\
           &     \tilde t                    & \longmapsto &     y = (P_1,...,P_k)(\tilde t) 
\end{array}                                      
\end{displaymath} 
L'équation de $\partial \Om_F$ en $z_0=\tilde \si \lbrace x_0,u_0 \rbrace$ devient :
\[ \partial \Om_F = \lbrace (y,w)  \in V_0  \times (\C,v_0) \, , \, \delta_0 . e ^ {- { \pi \over 2} H(\Phi ^{-1}(y),\Phi ^{-1}(y))} . \abs v = 1 \rbrace \]
Si on pose $w := {2 \over \pi} \left ( \log (\delta_0.v) - i \Im \log (\delta_0.v_0) \right)$, cette équation s'écrit :
\[ { \partial \Om_F = \lbrace (y,w) \in V_0 \times (\C,0) \, , \,  \Re (w) - H(\Phi ^{-1}(y),\Phi ^{-1}(y)) = 0  \rbrace } \]
Cela termine la description de $\partial \Om_F$.\\

On retrouve la description de \cite{BL1}. En effet, en dimension $1$, le stabilisateur d'un point critique de $\si$ est cyclique d'ordre $m=2,3,4$ ou $6$. L'application $\Phi$ est alors donnée par $\tilde t \mapsto y={\tilde t}^m$ et l'équation de $\partial \Om_F$ s'écrit $ \partial \Om_F = \lbrace (y,v) \in (\C,0) \times (\C,0), {\abs y}^{2 \over m} = \Re (v) \rbrace$.

\section{Exemple de Lattès en dimension $2$}

\subsection{Les couples $(A^2,G)$ tels que $A^2 / G \simeq \Pj^2$}\label{class}

La liste de ces couples et les valeurs critiques du passage au quotient $A^2 \to \Pj^2$ sont dus à Kaneko, Tokunaga et Yoshida (cf \cite{KTY}) : 
\begin{enumerate}
\item $\left(A_\om \times A_\om, G(2,1,2)\right)$, $4$ droites et une quadrique
\item $\left(A_\rho \times A_\rho, G(3,1,2)\right)$, $3$ droites et une quadrique
\item $\left(A_i \times A_i, G(4,1,2)\right)$, idem
\item $\left(A_\rho \times A_\rho, G(6,1,2)\right)$, idem
\item $\left(A_i \times A_i, \langle G(4,2,2),{{1+i}\over2} \binom{1}{1} \rangle\right)$, $6$ droites
\item $\left(A_\om \times A_\om, \SS_3\right)$, la courbe duale d'une cubique non singulière
\end{enumerate}
où $A_\om$ désigne le tore associé au réseau $\Z + \om\Z$, $G(m,p,2)$ le groupe engendré par les matrices : 
\[ \left(
\begin{array}{cc}
                   0       & 1 \\
                   1       & 0  
\end{array}
                \right)\  
\left(
\begin{array}{cc}
     0                     &   e ^{ {2i\pi \over m}} \\
     e ^{ -{2i\pi \over m}}  &   0
\end{array}
                \right)\  
\left(
\begin{array}{cc}
     e ^{ {2ip\pi \over m}} & 0 \\
             0             & 1 
\end{array}
                \right)\] 
et $\SS_3$ la représentation de $S_3$ suivante :
\[ \left(
\begin{array}{cc}
                   1       & 0 \\
                   0       & 1  
\end{array}
                \right)\ 
\left(
\begin{array}{cc}
                   -1       & -1 \\
                    0       & 1  
\end{array}
                \right)\  
\left(
\begin{array}{cc}
     0                     &   1 \\
     1                &   0
\end{array}
                \right)\ 
\left(
\begin{array}{cc}
     1                     &   0 \\
    -1                &   -1
\end{array}
                \right)\  
\left(
\begin{array}{cc}
            -1         & -1 \\
             1             & 0 
\end{array}
                \right)\ 
\left(
\begin{array}{cc}
            0         & 1 \\
             -1            & -1 
\end{array}
                \right)  
\]

Toutes ces situations induisent des exemples de Lattès (prendre $D=n.$Id). A titre d'exemple, donnons une équation de la singularité de $\partial \Om_F$ en $z_0 = \tilde \si \{ 0,u_0 \}$ pour le premier couple. D'après le théorème $3.13$ de \cite{F}, une base de $\C[X,Y]^G$ consiste en la donnée de deux polyn\^omes homogènes $(P,Q)$ dont la jacobienne n'est pas identiquement nulle, et tels que le produit de leur degré soit égal au cardinal du groupe $G$. Les polyn\^omes $P(X,Y)=X^2+Y^2$ et $Q(X,Y)=X^2 Y^2$ vérifient ces conditions. On obtient l'expression :
\begin{displaymath}
  \Phi :
\begin{array}{rccc}
           & (\C^2,0)/G             &  \to    &   (\C^2,0)           \\
           &     (X,Y)                     & \mapsto &   (\theta_1=X^2+Y^2,\theta_2=X^2 Y^2) 
\end{array}                                      
\end{displaymath} 
Les relations coefficients-racines d'un polyn\^ome du second degré fournissent $X^2,Y^2 \in \lbrace { - \theta_1 \pm \sqrt{\theta_1^2 - 4 \theta_2}  \over 2 } \rbrace$. La forme $H$ est un multiple de la forme standard, puisqu'elle est invariante par $G(2,1,2)$ (cf lemme \ref{herm}). Une équation de la singularité s'écrit donc 
\[  \partial \Om_F = \lbrace (\theta,v) \in (\C^2,0) \times (\C,0), \bigabs { \theta_1 + \sqrt{\theta_1^2 - 4 \theta_2} }   + \bigabs { \theta_1 - \sqrt{\theta_1^2 - 4 \theta_2} }  =\Re (v) \rbrace  \]

\subsection{Quelques exemples de Lattès}\label{UFS}

Un endomorphisme holomorphe $f$ de $\Pj^k$ est dit critiquement fini si l'orbite $\bigcup_{n \geq 1} f^n(\CC_f)$ de son ensemble critique $\CC_f$ est algébrique. Un exemple de Lattès est une application critiquement finie :
\begin{lem}\label{cf}
Soit $f$ un exemple de Lattès de $\Pj^k$ vérifiant $\si \circ D = f \circ \si$. Alors $PC(f)$ est contenu dans l'ensemble algébrique des valeurs critiques de $\si$.
\end{lem}

\proof
Pour tout $n$, on a $\si \circ D^n = f^n \circ \si$. Si $x \in \CC_f$, $f^n \circ \si$ branche aux points de $\si^{-1} \{ x \}$. Puisque $D^n$ est un rev\^etement, $\si$ branche aux points de $D^n \si^{-1} \{ x \}$ et $f^n(x)$ est une valeur critique de $\si$.\fin 

\begin{prop}
Les endomorphismes critiquement finis suivants sont des exemples de Lattès :
\begin{displaymath}
\begin{array}{rccc}
                  & \Pj^2  &  \longrightarrow    &  \Pj^2 \\
        f_1:       & [x:y:z]  & \longmapsto & [(-x+y+z)^2:(x-y+z)^2:(x+y-z)^2] \\       
        f_2:       & [x:y:z]  & \longmapsto & [(x-y+z)^2:(-x+y+z)^2:(x+y-z)^2] \\
        f_3:       & [x:y:z]  & \longmapsto & [(x+y-z)^2:(-x+y+z)^2:(x-y+z)^2] 
\end{array}                                      
\end{displaymath}
Ils sont semi-conjugués par la situation $5$ de la classification aux dilatations :
\[ D_1 = \sqrt2 \left(
\begin{array}{cc}
     e ^{ {i\pi \over 4}}  & 0 \\
                   0       & e^{ {i\pi \over 4}}  
\end{array}
                \right),\ 
 D_2 = \sqrt2 \left(
\begin{array}{cc}
    {1 \over \sqrt 2} &   {1 \over \sqrt 2} \\
    {1 \over \sqrt 2} &  {-1 \over \sqrt 2}
\end{array}
                \right),\ 
D_3 = \sqrt2 \left(
\begin{array}{cc}
    {1 \over \sqrt 2} &   {1 \over \sqrt 2} \\
    {i \over \sqrt 2} &  {-i \over \sqrt 2}
\end{array}
                \right)\]  
\end{prop}
\ 

\proof
On note $\wp$ la fonction de Weierstrass associée au réseau $\Z + i \Z$ et $\theta [ a \, b ]$ les fonctions th\^eta de Riemann pour ce m\^eme réseau. Leurs diviseurs sont respectivement $2[{1+i \over 2}] -2 [0]$ et $[ i(a+{1 \over 2}) + (b+ {1 \over 2}) ]$ (cf \cite{D}, Chap.II, §$2$). Soit $\theta_{jk}$ la fonction $\theta [ {j \over 2} \ {k \over 2} ]$. Puisque $\theta_{00}^2$ et $\theta_{11}^2$ sont de m\^eme type, il existe $c \in \C$ tel que $\wp = c.\theta_{00}^2.\theta_{11}^{-2}$. L'application de passage au quotient de la situation $5$ s'écrit :
\begin{displaymath}
\begin{array}{c}
     \si (x,y) =  [\left( c^2.\theta_{00}^2(x).\theta_{00}^2(y) + \al^2.\theta_{11}^2(x).\theta_{11}^2(y) \right)^2 :  \left( c^2.\theta_{00}^4(x)-\al^2.\theta_{11}^4(x)\right). \\
                                         
          \hspace {4 cm} \left( c^2.\theta_{00}^4(y)-\al^2.\theta_{11}^4(y) \right) : \left( c^2.\theta_{00}^2(x).\theta_{00}^2(y) -  \al^2.\theta_{11}^2(x).\theta_{11}^2(y) \right)^2  ]
\end{array}                                      
\end{displaymath} 
C'est l'application donnée dans \cite{KTY}, exprimée à l'aide de fonctions th\^eta, et composée par l'automorphisme $[X:Y:Z] \mapsto [X:X-\al^2 Y:Z]$ de $\Pj^2$. La démonstration consiste à vérifier les relations $\si \circ D_i = f_i \circ \si$. Nous aurons besoin des identités suivantes :
\begin{enumerate}
\item  $ \wp_x.\wp_{x + ix} =-{i \over 2} . \left( \wp_{x} ^ 2  - \al ^ 2 \right)$
\item  $( \wp_x - \wp_y )^2 .\wp_{x+y} . \wp_{x-y} = \left( \wp_x. \wp_y + \al^2 \right)^2$
\item  $( \wp_x - \wp_y )^2 .\left( \wp_{x+y} + \wp_{x-y} \right) = 2. \left( \wp_x + \wp_y \right)  \left( \wp_x.\wp_y - \al^2 \right)$
\end{enumerate}
Elles proviennent des relations ${\wp_x'} ^2 = 4.\wp_x. (\wp_x ^2 - \al ^2)$, $\wp_{ix}=-\wp_x$, $\wp_{-x}=\wp_x$, et de la formule d'addition :
   \[ ( \wp_x - \wp_y )^2.\wp_{x + y} = {1 \over 4} ( \wp'_x - \wp'_y )^2  - (\wp_x + \wp_y) ( \wp_x - \wp_y )^2 \]  
Vérifions la première. On a : 

\noindent $\wp_x.\wp_{x + ix}.( \wp_x - \wp_{ix} )^2 = \wp_x.{1 \over 4} ( \wp'_x - \wp'_{ix} )^2  - \wp_x.(\wp_x + \wp_{ix}) ( \wp_x - \wp_{ix} )^2 = \wp_x . { (1-i)^2 \over 4 } . {\wp_x'}^2 = $

\noindent $-2i.{\wp_x}^2.(\wp_x^2 - \al^2)$. 
La conclusion vient de l'autre identité $\wp_x.\wp_{x + ix}.( \wp_x - \wp_{ix} )^2 = \wp_x.\wp_{x + ix}.4.\wp_x ^2$.\\ 

\noindent Passons à la deuxième : 

\noindent  $( \wp_x - \wp_y )^2 .\wp_{x+y} . ( \wp_x - \wp_y )^2 . \wp_{x-y} = ( \wp_x - \wp_y )^2 .\wp_{x+y} . ( \wp_x - \wp_{-y} )^2 . \wp_{x-y} =$

\noindent  $\left ( {1 \over 4} (\wp'_x - \wp'_y)^2 -  (\wp_x + \wp_y) (\wp_x - \wp_y)^2  \right) \left( {1 \over 4} (\wp'_x - \wp'_{-y})^2 -  (\wp_x + \wp_{-y}) (\wp_x - \wp_{-y})^2  \right) = $

\noindent $\left ( {1 \over 4} ({\wp'_x}^2 + {\wp'_y}^2) -  (\wp_x + \wp_y) (\wp_x - \wp_y)^2  \right)^2 - \left( {1 \over 2} . \wp'_x . \wp'_y \right)^2 = $ 

\noindent $\left(  \wp_x. (\wp_x ^2 - \al ^2) + \wp_y. (\wp_y ^2 - \al ^2) - (\wp_x + \wp_y) (\wp_x - \wp_y)^2  \right)^2  - 4 \wp_x . \wp_y.(\wp_x ^2 - \al ^2)(\wp_y ^2 - \al ^2) = $

\noindent $\left(  \wp_x^3 - \al ^2.\wp_x + \wp_y ^3 - \al ^2.\wp_y - (\wp_x + \wp_y)(\wp_x^2 - 2 \wp_x.\wp_y + \wp_y^2) \right)^2  $

\noindent $ \hspace{ 9.6 cm } - 4 \wp_x . \wp_y.(\wp_x ^2 - \al ^2)(\wp_y ^2 - \al ^2) = $ 

\noindent $\left( \wp_x.(\wp_y ^2 - \al ^2) + \wp_y.(\wp_x ^2 - \al ^2) \right) ^2  - 4 \wp_x . \wp_y.(\wp_x ^2 - \al ^2)(\wp_y ^2 - \al ^2) = $

\noindent $  \left( \wp_x.(\wp_y ^2 - \al ^2) - \wp_y.(\wp_x ^2 - \al ^2) \right) ^2 = \left( \wp_x - \wp_y \right)^2 \left( \wp_x. \wp_y + \al^2 \right)^2 $\\

\noindent La dernière relation provient de :

\noindent $ ( \wp_x - \wp_y )^2 .\left( \wp_{x+y} + \wp_{x-y} \right) = ( \wp_x - \wp_y )^2 . \wp_{x+y} + ( \wp_x - \wp_{-y} )^2.\wp_{x-y}  = $

\noindent $ {1 \over 4} (\wp'_x - \wp'_y)^2 - (\wp_x + \wp_y).( \wp_x - \wp_y )^2 + {1 \over 4} (\wp'_x - \wp'_{-y})^2 - (\wp_x + \wp_{-y}).( \wp_x - \wp_{-y} )^2 = $

\noindent $ 2. \left ( {1 \over 4}. ( {\wp'}^2_x + {\wp'}^2_y ) - (\wp_x + \wp_y).( \wp_x - \wp_y )^2 \right) = $

\noindent $ 2. \left(  \wp_x. (\wp_x ^2 - \al ^2) + \wp_y. (\wp_y ^2 - \al ^2) - (\wp_x + \wp_y) (\wp_x - \wp_y)^2  \right) $

\noindent On retrouve cette expression dans les calculs de la relation $2$. Elle est égale à :  

\noindent $ 2. \left( \wp_x.(\wp_y ^2 - \al ^2) + \wp_y.(\wp_x ^2 - \al ^2) \right) ^2 =  2. \left( \wp_x + \wp_y \right)  \left( \wp_x.\wp_y - \al^2 \right)$\\

\textbf{Cas de l'endomorphisme $f_1$} :

\noindent Le diviseur des fonctions th\^eta $c.\theta_{00}^2((1+i).x)$ et $c^2. \theta_{00}^4(x) - \al^2.\theta_{11}^4(x)$ vaut $2 [ {1 \over 2} ] + 2 [ {i \over 2} ]$. Celui de $\theta_{11}^2((1+i)x)$ et $c.\theta_{00}^2(x).\theta_{11}^2(x)$ vaut $2[0] + 2[ {1+i \over 2}]$. Il existe donc des fonctions th\^eta $\varphi_j$ telles que :
\begin{displaymath}
\begin{array}{rcl}
    c.\theta_{00}^2((1+i)x) & =  & \varphi_0(x).\left( c^2. \theta_{00}^4(x) - \al^2.\theta_{11}^4(x) \right) \\  
     \theta_{11}^2((1+i)x)   & = &  \varphi_1(x).c.\theta_{00}^2(x).\theta_{11}^2(x)
\end{array}                                      
\end{displaymath}  
En divisant ces deux expressions, la relation $1$ donne $\varphi_0 = -{ i \over 2 } . \varphi_1$. Les premières et dernières coordonnées de $\si \circ D_1$ s'écrivent alors : ($\epsilon$ vaut respectivement $-1,1$)
\begin{displaymath}
\begin{array}{l}
     \left( -{1 \over 4}.\varphi_1(x).\varphi_1(y) \right)^2   \lbrace (c^2. \theta_{00}^4(x) - \al^2.\theta_{11}^4(x))(c^2. \theta_{00}^4(y) - \al^2.\theta_{11}^4(y)) \\
                                \\
  \hspace {8 cm}    4.\epsilon.\al^2.c^2.\theta_{00}^2(x).\theta_{11}^2(x).\theta_{00}^2(y).\theta_{11}^2(y) \rbrace^2 \\
\end{array}                                      
\end{displaymath}   
et coïncident avec $ \left( -{1 \over 4}.\varphi_1(x).\varphi_1(y) \right)^2 . (\si_2 + \epsilon.(\si_3 - \si_1))^2$, où $\si_k$ est la $k$-ième coordonnée de $\si$. La deuxième coordonnée est égale à :
\begin{displaymath}
\begin{array}{l}
     \left( -{1 \over 4}.\varphi_1(x).\varphi_1(y) \right)^2 .  \lbrace ( c^2.\theta_{00}^4(x) - \al^2.\theta_{11}^4(x))^2 + 4.\al^2.\left( c.\theta_{00}^2(x).\theta_{11}^2(x) \right)^2 \rbrace. \\
  \hspace {5 cm}  \lbrace ( c^2.\theta_{00}^4(y) - \al^2.\theta_{11}^4(y))^2 + 4.\al^2.\left( c.\theta_{00}^2(y).\theta_{11}^2(y) \right)^2 \rbrace = \\
                                        \\
  \hspace {1 cm}   \left( -{1 \over 4}.\varphi_1(x).\varphi_1(y) \right)^2 . \left(  \left( c^2.\theta_{00}^4(x) + \al^2.\theta_{11}^4(x) \right) . \left(  c^2.\theta_{00}^4(y) + \al^2.\theta_{11}^4(y) \right) \right)^2
\end{array}                                      
\end{displaymath} 
et coïncide avec $ \left( -{1 \over 4}.\varphi_1(x).\varphi_1(y) \right)^2 . (\si_1 - \si_2 + \si_3)^2$.\\

\textbf{Cas de l'endomorphisme $f_2$} :

\noindent On écrit $\si$ sous la forme :
\[ \si (x,y) = [ ( \wp_x.\wp_y + \al ^2 )^2 : ( \wp^2 _x - \al^2 ) ( \wp^2 _y - \al^2 ) :  ( \wp_x.\wp_y - \al ^2 )^2 ] \]   
Avec la relation $2$, les premières et dernières coordonnées de $D_2 \circ \si$ multipliées par $\left( \wp_x - \wp_y \right)^4$ s'écrivent : ($\epsilon$ vaut respectivement $1,-1$) \[ \left( \wp_x - \wp_y \right)^4     \left( \wp_{x+y} . \wp_{x-y} +\epsilon.\al^2 \right)^2 = \left( \left( \wp_x. \wp_y + \al^2 \right)^2 + \epsilon.\al^2.\left( \wp_x - \wp_y \right)^2 \right)^2 \] 
et coïncident avec $\left( \si_1 + \epsilon.( \si_3 - \si_2 ) \right)^2$. On utilise ensuite les relations $2$ et $3$ pour la deuxième coordonnée multipliée par $\left( \wp_x - \wp_y \right)^4$ :
\begin{displaymath}
\begin{array}{l} 
 \left( \wp_x - \wp_y \right)^4.(\wp_{x+y} - \al).(\wp_{x-y} - \al).(\wp_{x+y} + \al) . (\wp_{x-y} + \al) = \\
 
\left( \wp_x - \wp_y \right)^4.\left(  (\wp_{x+y}.\wp_{x-y} + \al^2) - \al. (\wp_{x+y} + \wp_{x-y}) \right).\\

 \hspace{5 cm}\left(  (\wp_{x+y}.\wp_{x-y} + \al^2) + \al. (\wp_{x+y} + \wp_{x-y}) \right) = \\
 
\left( \wp_x - \wp_y \right)^4.\left(  \left( \wp_{x+y}.\wp_{x-y} + \al^2 \right) ^2 - \al^2 . \left( \wp_{x+y} + \wp_{x-y} \right)^2 \right) = \\
 \left(   (\wp_x.\wp_y + \al^2)^2  +  \al^2.( \wp_x - \wp_y )^2  \right)^2  - 4 \al^2 .  \left( \wp_x + \wp_y \right)^2  \left( \wp_x.\wp_y - \al^2 \right)^2 = \\
 \left(   (\wp_x.\wp_y - \al^2)^2  +  \al^2.( \wp_x + \wp_y )^2  \right)^2  - 4 \al^2 .  \left( \wp_x + \wp_y \right)^2  \left( \wp_x.\wp_y - \al^2 \right)^2 = \\
 \left(   (\wp_x.\wp_y - \al^2)^2  -  \al^2.( \wp_x + \wp_y )^2  \right)^2  = (\si_3 + (\si_2 - \si_1))^2
\end{array}                                      
\end{displaymath}

\ 

\textbf{Cas de l'endomorphisme $f_3$} :

\noindent On le déduit de l'expression de $f_2$ et de l'identité $\wp_{i(x-y)} = - \wp_{x-y}$. Les premières et dernières coordonnées sont permutées, la deuxième reste inchangée.\fin

\begin{rem}
\rm Des calculs analogues sont menés dans \cite{U1} pour l'endomorphisme $f_1$. Notre apport est de traiter les endomorphismes $f_2$, $f_3$ et d'utiliser \cite{KTY} afin de montrer que $f_1$, $f_2$ et $f_3$ sont des exemples de Lattès.
\end{rem}

\subsection{Un endomorphisme critiquement fini qui n'est pas un exemple de Lattès}

\begin{prop} 
L'endomorphisme critiquement fini :
\begin{displaymath}
g :
\begin{array}{rccc}
                  & \Pj^2  &  \to    &  \Pj^2 \\
                  & [x:y:z]  & \mapsto & [(x-2y)^2:(x-2z)^2:x^2]        
      
\end{array}                                      
\end{displaymath}
n'est pas un exemple de Lattès.
\end{prop}

\proof
Regardons la dynamique de l'ensemble critique de $g$ dont les composantes sont $\{ X=0 \}, \{ 2Y=X \}, \{ 2Z=X \}$ : 
\[
\xymatrix{
    \{ X=0 \}  \ar[r]  &  \{ Z=0 \} \ar[r]  &  \{ Y=Z \} \ar[r]  \ar[d]  &  \{ X=Y \} \\
    \{ 2Y=X \} \ar[u]  &                    &  \{ X=Z \} \ar[ru]         &            \\
    \{ 2Z=X \} \ar[r]  &  \{ Y=0 \} \ar[ru] &                            &           
          }
\]
L'orbite post critique contient $6$ droites. Si c'est un exemple de Lattès, il est nécessairement associé à la cinquième situation de la classification du paragraphe \ref{class}, d'après le lemme \ref{cf}. Supposons alors que $g$ fasse commuter :
\[
\xymatrix{
          A_i \times A_i \ar[r]^{D} \ar[d]_\si     & A_i \times A_i  \ar[d]^\si \\
          \Pj^2 \ar[r]^{g}                         & \Pj^2
        }
\]
On a donc modulo $\langle G(4,2,2),{{1+i}\over2} \binom{1}{1} \rangle$ : 
\[
\xymatrix{
          \lbrace (x,x+{1+i \over 2}) \rbrace\ar[r]^{D} \ar[d]_\si     & \lbrace (x,ix+{1+i \over 2}) \rbrace\ar[d]^\si \\
          \{ X=0 \}  \ar[r]^{g}                                        & \{Z=0\}
        }
\]
Le stabilisateur de $\lbrace (x,ix+{1+i \over 2}) \rbrace$ est d'ordre $2$, engendré par  $\left(
\begin{array}{cc}
     0  &   -i  \\
     i  &    0
\end{array}
                \right)$. On obtient alors une contradiction lorsque l'on compare les multiplicités dans le diagramme précédent : $\si \circ D$ est génériquement d'ordre $2$ au voisinage de $\lbrace (x,x+{1+i \over 2}) \rbrace$, alors que $g \circ \si$ est d'ordre supérieur, car $ \{ X=0 \} $ est dans l'ensemble critique de $g$. \fin

\ 

Christophe Dupont

Laboratoire Emile Picard, UMR 5580

Université Toulouse III

118 Route de Narbonne, 31062 Toulouse Cedex 04, France

dupont@picard.ups-tlse.fr

Fax : (33) 05 61 55 82 00
\end{document}